\documentclass[11pt]{amsart}
\usepackage{amsmath,amsthm,latexsym,amssymb}
\usepackage[dvips]{graphicx}
\usepackage{color}
\definecolor{bleu}{rgb}{0.00,0.4,0.90}
\definecolor{magenta}{rgb}{1.0, 0.0, 1.0}
\newcommand{\nc}{\normalcolor}

\numberwithin{equation}{section}

\newtheorem{theorem}{Theorem}

\newtheorem{remark}{Remark}
\newtheorem{definition}{Definition}
\numberwithin{theorem}{section}

\numberwithin{corollary}{section}
\numberwithin{lemma}{section}
\numberwithin{definition}{section}
\numberwithin{proposition}{section}
\numberwithin{remark}{section}

\newcommand{\R}{\mathbb R}

\newcommand{\medint}{-\kern  -,375cm\int}

\newcommand{\be}{\begin{equation}}
\newcommand{\ee}{\end{equation}}
\newcommand{\beqn}{\begin{eqnarray}}
\newcommand{\eeqn}{\end{eqnarray}}
\newcommand{\de}{\partial}

\begin{document}

\title[On the behavior of the solutions]{Uniqueness, regularity and behavior in time of the solutions to nonlinear  anisotropic parabolic equations }


\author[Di Blasio]{Giuseppina Di Blasio}
\address{Universit\`{a} degli Studi della Campania \\
Dipartimento di Matematica e Fisica, \textquotedblleft L. Vanvitelli\textquotedblright , \\
Viale Lincoln, 5 - 81100
Caserta, Italy}
\email{giuseppina.diblasio@unicampania.it}

\author[Porzio]{Maria Michaela Porzio}
\address{Sapienza Universit\`a di Roma\\
Dipartimento di Pianificazione, Design, Tecnologie del\-l'Ar\-chi\-tet\-tu\-ra\\
via Flaminia 72, 00196 Roma, Italy}
\email{mariamichaela.porzio@uniroma1.it}

\maketitle

\begin{abstract}
 In this paper  we study a class of nonlinear anisotropic parabolic problems in  bounded domains. In detail, we study the influences of the initial data and the forcing term $f$ on the behavior of the solutions. We prove existence and uniqueness results. We indagate on the behavior in time of the solutions with a particular attention to the autonomous case $f(x,t)=f(x)$.
\end{abstract}


\section{Introduction}

In this paper we study a class of nonlinear anisotropic parabolic
problems whose model case is

\begin{equation}\label{anis}
 \left\{\begin{array}{ll}
u_t- \sum_{i=1}^N \partial_{x_{i}}\left( |u_{x_i}|^{p_i-2}u_{x_i} \right)  = f & {\rm in} \quad \Omega_{T}:=\Omega\times\left(  0,T\right) ,\\
u=0 & {\rm on} \quad \partial\Omega \times
(0,T),\\
u(x,0)=u_0(x) & {\rm on} \quad \Omega\,,
\end{array}
\right.
\end{equation}
where $\Omega$ is an open, bounded subset of $\mathbb{R}^{N}$, $N\geq2,$ $T>0$ and $p_i>1$, for every $i=1,...,N$,   such that their harmonic mean
$\bar{p}$ belongs to $(1,N)$ and the data $f$ and $u_0$ are measurable functions.

In the case $p_i=p $ for every $i=1,...,N$ the equation in \eqref{anis} is known as the parabolic orthotropic $p$-Laplacian equation and even if its principal part is homogeneous, it fails to be isotropic. It becames the isotropic $p$-Laplacian only in the case $p_i=2$ for every $i=1,...,N$.

In the last decades the studies of anisotropic problems represent a very active research area. The increasing interest in this kind of matter is justified by their applications in many areas. Indeed many materials, such as liquid crystals, wood or earth’s crust, usually present different
diffusion rates along different directions. So problems \eqref{anis} provides the mathematical models for many natural phenomena.
For example, they are the mathematical description of the dynamics of fluids in anisotropic media when the conductivities of the
media are different in different directions or in biology they describes a model
for the propagation of epidemic diseases in heterogeneous domains (see, for example, \cite{ADS}, \cite{BK}).

In the classical theory of regularity for solutions to anisotropic elliptic problems the anisotropy
condition depends on differential operators whose growth with respect to the
partial derivatives of $u$ is governed by different powers (see \emph{e.g.}
\cite{antontsev-chipot-08,BMS, BaCri,Barbara, dBFZ, dBFZ_2, castro, ELM, FGK,FGL,Gi, LI,Mar, MRSC,Str} and references therein).

The anisotropy prescribed by \eqref{anis} falls into the wider class of problems with nonstandard
growth condition. There are many papers related to problems governed by fully anisotropic growth conditions and also by
different type of anisotropy (see \emph{e.g.} \cite{AFTL, BFK, cianchi immersione, cianchi
anisotropo, dBL} and references therein).

Regarding the parabolic anisotropic problems, much less is known. Several results are contained for example in \cite{AdBF},\cite{BK0}\cite{DMV},\cite{FVV}, \cite{Li}, \cite{LiZ}, \cite{M}, \cite{P3}, but the results on  the solutions of \eqref{anis} are far to be complete and many interesting questions are open.

Aim of this paper is to answer to some open problems concerning this kind of problems to reach  a deeper knowledge of the behavior of   the solutions of \eqref{anis}.

For example, it is well known that if the data $u_0$ and $f$ are summable functions and the coefficients $p_i$ satisfy
$$
2 -\frac{1}{N+1} \, < \, p_i \, < \, \overline{p} \, \frac{N+1}{N}\,, \quad  \quad \overline{p} < N,
$$
then there exists at least a weak solution of \eqref{anis} (in the sense of Definition \ref{defsol} below) (see \cite{LiZ}). Anyway  no uniqueness results seems to be known on our problem \eqref{anis} and it is worth to be noticed that
 if the data $u_0$ and $f$ are only summable functions we cannot expect that the weak solutions of \eqref{anis} are unique because this lack of uniqueness still appears in the particular case of the othotropic $p$-Laplacian equation together with other similar parabolic problems like  the classical isotropic $p$-Laplacian equation.

 To guarantee the uniqueness of solutions, for example in these last two cases of othotropic $p$-Laplacian equation or of the $p$-Laplacian equation (but also for other type of parabolic equations) further requirements on the solutions are needed in order to choose ``the only one solution''  which has this further property. Hence, different type of requirements on the solutions have been proposed in literature to guarantee the uniqueness of solutions like, for example, the entropy solutions, the renormalized solutions or the solutions obtained by approximation. Here we have preferred to study the solutions obtained by approximations (see Definition \ref{defappro} below) since they are the more natural  way to construct solutions.

Hence our first goal has been to understand if also in the anisotropic case the further requirement to be a solution ``constructed by approximation'' guarantees its  uniqueness.

We point out that in the case of the othotropic $p$-Laplacian equation it is known that there exists and it is unique the solution $u$ obtained by approximation (see \cite{debmon}).

In Theorem \ref{thuni} we prove that there exists one and only one solution $u \in C([0,T];L^1(\Omega))$ constructed by approximation of \eqref{anis}. Moreover, there is also a sort of continuous dependence from the data, since if $u$ is the unique solution in $C([0,T];L^1(\Omega))$ constructed by approximation of \eqref{anis}, for every $t \in (0,T)$ it results
\begin{equation*}\label{cont}
\|u(t)-v(t)\|_{ L^1(\Omega)} \leq \|u_0-v_0\|_{ L^1(\Omega)} + \|f-g\|_{ L^1(\Omega \times (0,t))}\,,
\end{equation*}
where  $v $ is the unique solution in  $C([0,T];L^1(\Omega))$ constructed by approximation of the following problem
\begin{equation}\label{manisv}
 \left\{\begin{array}{ll}
v_t- \sum_{i=1}^N \partial_{x_{i}}\left( |v_{x_i}|^{p_i-2}u_{x_i} \right)  = g & {\rm in} \quad \Omega_{T}\,,\\
v=0 & {\rm on} \quad \partial\Omega \times
(0,T),\\
v(x,0)=v_0(x) & {\rm on} \quad \Omega\,,
\end{array}
\right.
\end{equation}
with $v_0 \in L^1(\Omega)$ and $g \in L^1(\Omega_T)$.

Moreover, this result can be extended also to the so called global solutions (see Definition \ref{global}).

Once understood that our problem admits a unique solution obtained by approximation, we have studied the regularity of this solution and its behavior in time.

We point out that there are several results in literature on the regularity  of the solution of \eqref{anis} in the sense of summability in Lebesgue space. For example in \cite{P3} the author proves that if
\begin{equation}\label{ASani}
f \in  L^m(\Omega_T)\,, \quad m > 1 + \frac{N}{\overline{p}}\,, \quad \hbox{\rm and} \quad  u_0 \in  L^{\infty}(\Omega)
\end{equation}
then the solution $u$ is bounded in $\Omega_T$. A natural question is to understand what happens if the forcing term $f$ satisfies \eqref{ASani}  but $u_0 \not\in  L^{\infty}(\Omega)$.

In the isotropic case (see \cite{decay}) there is a partial answer to this question when the datum $f \equiv 0$ and  $u_0 $ is only a summable function. In this case the solution becomes immediately bounded and satisfies
 the following bound
\begin{equation}\label{decayf0}
\|u(t)\|_{L^{\infty}(\Omega)} \leq C \, \frac{\|u_0\|^{\frac{\overline{p}}{N(\overline{p}-2)+ \overline{p} }}_{L^{1}(\Omega)}}{t^{\frac{N}{N(\overline{ p}-2)+ \overline{p} }}} \,.
\end{equation}

We recall that estimates like \eqref{decayf0} are referred in literature as decay estimates or ultracontractive estimates and
 are known to appear for many different type of evolution equations like the $p$-Laplacian equation, the porous media equations, the fast diffusion equations, the doubly nonlinear equation together with other type of degenerate or singular parabolic problems (see for example  \cite{AP}, \cite{BCP}, \cite{BG}, \cite{CG},  \cite{DiBH1}, \cite{DiBH2}, \cite{FVV},  \cite{HV}, \cite{decay},   \cite{V1}  and the references therein). The interest to this kind of estimates is due to the fact that not only they imply that, even if the initial datum  $u_0$ is unbounded then the solution $u(t)$ becomes immediately bounded, but they describe the behavior in time of the solutions revealing that they decay to zero in the $L^{\infty}$-norm as $t$ goes to infinity. Indeed, this kind of bounds describes also the blow-up of the $L^{\infty}$-norm of the solution $u(t)$ as soon as $t$ goes to zero if the initial data are unbounded.

In the anisotropic case we prove that a similar behaviour in time occurs. More precisely, we show  (see Theorem \ref{asinto}) that if  $f$ satisfies \eqref{ASani} with $f \not\equiv 0$ and  $u_0$ is not bounded and not regular, that is $u_0 \in  L^{1}(\Omega)$, then again as soon as $t>0$ the solution $u(t)$ of \eqref{anis} becomes immediately bounded. It appears that the immediate boundedness seems to be produced by the regularity of the forcing term $f$ and not by the regularity of the initial datum $u_0$.

 In the paper we analyze deepest this interesting property of the solutions of the anisotropic problem \eqref{anis} and we get that as soon as $t>0$ the regularity of initial datum $u_0$ doesn't influence the regularity of the solution $u$ of \eqref{anis} since such a solution has the same summability regularity of  the solution $v$ of the same problem but with a bounded initial datum (see Theorem \ref{regLrs}).

Finally, we study the asymptotic behavior of the solutions. We  prove that if $u$ and $v$ are the unique global solutions constructed by approximation of, respectively, \eqref{anis} and \eqref{manisv} then they have the same behavior for $t$ large. In detail, the difference $u-v$ belongs to $L^{\infty}(\Omega \times [t,+\infty]) $ (for every $t>0$) and it results
$$
\lim_{t \rightarrow +\infty} \|u-v\|_{L^{\infty}(\Omega \times [t,+\infty])} = 0\,.
$$
We point out that the previous results are rather surprising since it affirms that the difference $u-v$ belongs to $L^{\infty}(\Omega \times [t,+\infty])$   even if both the solutions $u$ and $v$ are generally unbounded (see Theorem \ref{asinto}).

  We consider also the  autonomous case $f(x,t) \equiv f(x)$ and we prove that in this case all the global solutions obtained by approximation of \eqref{anis}, whenever is the initial datum $u_0$ that they assume, for $t \rightarrow +\infty$ converge in the $L^{\infty}$-norm to the unique solution $w$ constructed by approximation of the associated stationary problem
\begin{equation}\label{ellani_intro}
\left\{
\begin{array}
[c]{lll}%
\sum_{i=1}^N \partial_{x_{i}}\left( |w_{x_i}|^{p_i-2}w_{x_i} \right) = f(x) & \quad \hbox{\rm in }%
\Omega\,,  & \\
 w(x)=0 & \quad \hbox{\rm on }\;\partial\Omega\,. &
\end{array}
\right.
\end{equation}

We point out that also for the anisotropic problem \eqref{ellani_intro} the solutions are not unique if the forcing term $f$ is only assumed to be a summable function. Hence we have shown that again the requirement that the solution $w$ of \eqref{ellani_intro} is constructed by approximation guarantees that it is also unique. Notice that also in this case the convergence of $u$ to $w$ in the $L^{\infty}$-norm is rather surprising since when the data $u_0$ and $f$ are only summable functions both  the solutions $v$ and $w$ are generally unbounded.

The paper is organized as follows: in Section \ref{secstate} we state all our results, in Section \ref{secpre}  we recall some useful results that we will use  in  Section \ref{secproofs} in  proving our results.


\section{Statement of results}\label{secstate}

In this paper we deal with a class of nonlinear parabolic problems subject to general growth
conditions and having the form
\begin{equation}\label{P_parab}
\left\{
\begin{array}
[c]{lll}%
u_{t}+\sum_{i=1}^N  \partial_{x_{i}}(a_i(x,t,\nabla u))=f(x,t) & \quad\text{in }%
\Omega_{T}  & \\
 u(x,t)=0 & \quad\text{on}\;\partial\Omega\times\left(  0,T\right)  , &\\
u(x,0)=u_{0}(x) & \quad\text{in}\;\Omega, &
\end{array}
\right.  %
\end{equation}
where $\Omega$ is an open, bounded subset of $\mathbb{R}^{N}$, $N\geq2,$ $T>0$.

Here,  $a_i:\Omega_{T}\times\mathbb{R}%
^{N}\rightarrow\mathbb{R}^{N}$, for every $i=1,...,N$, is a Carath\'{e}odory function  satisfying the following structure conditions

\begin{eqnarray}
\label{ii1}&& \alpha  |\xi_i|^{p_i}\leq   a_i(x,t, \xi )\, \xi_i \,, \quad \,\, p_i > 1 \,,  \\
\label{ii2}&& |a_i (x,t, \xi ) | \leq \beta_i  | \xi_i|^{p_i-1}  + h(x,t)^{1-\frac{1}{p_i}} \,, \\
\label{ii3}&& ( a_i(x,t, \xi ) - a_i(x,t, \eta ) )\, (\xi_i - \eta_i ) > 0\,, \quad \xi_i \not= \eta_i \,,
\end{eqnarray}

\noindent for a.e. $(x,t) \ \in \Omega_T$ and for any vector $\xi,\eta$  in $\mathbb R^N$, where $\alpha$ and $\beta_i$  are positive constants satisfying  $\alpha \leq \beta_i$   and $h \in L^1(\Omega_T)$.\nc

Moreover, we assume that

\begin{equation*}\label{dato_u_0}
u_0\in  L^{1}(\Omega)\quad \quad \text{and}\quad \quad f\in  L^{1}(\Omega_T).
\end{equation*}

\bigskip

To guarantee the existence of solutions to problem \eqref{P_parab} under the previous weak summability conditions on the data, the following assumption on the coefficients $p_i$ will be retained
\begin{equation}\label{assupi}
2 -\frac{1}{N+1} \, < \, p_i \, < \, \overline{p} \, \frac{N+1}{N}\,, \quad  \quad \overline{p} < N\,,
\end{equation}
 where  $\bar{p}$  is the harmonic mean of $\vec{p}=
(p_1,\cdots,p_N)$, \textit{i.e.}
\begin{equation}\label{p bar}
\frac{1}{\overline{p}}=\frac{1}{N}\sum_{i=1}^N\frac{1}{p_i}.
\end{equation}

Notice that in the isotropic case $p_i \equiv p$ the previous condition (see \cite{BG89}) becomes the classical assumption
$$
2 -\frac{1}{N+1} \, < \, p \, < \, N\,.
$$

\begin{remark}\label{rempbar}
We point out that the assumption $p_i > 2- \frac{1}{N+1}$  implies that  $ \bar p > 2- \frac{1}{N+1}$. Consequently, being $2- \frac{1}{N+1} > \frac{2N}{N+1}$, assumption
 \eqref{assupi} implies that    $ \bar p > \frac{2N}{N+1}$.
\end{remark}

\bigskip

%
%
%

We stress that hereafter, for the partial derivative $\partial_{x_i}u$ of a function $u$ we adopt the alternative
notation $u_{xi}$ .

Before stating our results we recall the definitions of solutions considered here.

\medskip

\begin{definition}\label{defsol}
We say that a function $u\in L^1(0,T;W^{1,1}_0(\Omega))$ is a weak solution of problem \eqref{P_parab} if $a_i(x,t,\nabla u)\in L^1(\Omega_T)$ for every $i=1,...,N$, and
\begin{equation*}
  \int\int_{\Omega_T}\left[-u \varphi_t + \sum_{i=1}^{N} a_i(x,t,\nabla u) \varphi_{ x_i} \right] \,dx dt=
  \int_{\Omega} u_0\varphi(0) \,dx+ \int\int_{\Omega_T}f \varphi \,dx dt
\end{equation*}
for every $\varphi\in C^{\infty}(\overline{\Omega_T})$ which is zero in a neighborhood of $(\partial\Omega \times
(0,T))\cup (\Omega \times \{T\})$.
\end{definition}

\begin{definition}\label{defappro}
We say that $u$ is a weak solution of  \eqref{P_parab} obtained by approximation (SOLA) if it is the a.e. limit in $\Omega_T$ of the solutions $u_n\in L^{\infty}(\Omega_T)\cap C([0,T];L^2(\Omega))\cap L^{\overrightarrow{p}}\left(  0,T;W_{0}%
^{1,\overrightarrow{p}}(\Omega)\right)$ of the following (approximating) problems
\begin{equation}\label{P_approx}
\left\{
\begin{array}
[c]{lll}%
(u_n)_{t}+\sum_{i=1}^N  \partial_{x_{i}}(a_i(x,t,\nabla u_n))=f_n(x,t) & \quad \hbox{\rm in }%
\Omega_{T},  & \\
 u_n=0 & \quad \hbox{\rm on }\;\partial\Omega\times\left(  0,T\right)  , &\\
u_n(x,0)=u_{0,n}(x) & \quad \hbox{\rm on }\;\Omega, &
\end{array}
\right.  %
\end{equation}
where $f_n$ and $u_{0,n}$ satisfy
\begin{equation*}
f_n\in L^{\infty}(\Omega_T), \quad f_n\rightarrow f \quad \text{in} \quad L^1(\Omega_T),
\end{equation*}
\begin{equation*}
u_{0,n}\in L^{\infty}(\Omega), \quad u_{0,n}\rightarrow u_0 \quad \text{in} \quad L^1(\Omega),
\end{equation*}
and  $L^{\overrightarrow{p}}\left(0,T;W_{0}^{1,\overrightarrow
{p}}(\Omega)\right)$ is  the classical anisotropic space whose  definition can be found in Section \ref{aniSpace}.
\end{definition}

\begin{definition}\label{global}
We say that $u$ is a global solution of \eqref{P_parab}, or that is the same, of a global solution of the following problem
\begin{equation}\label{gloanis}
 \left\{\begin{array}{ll}
u_t- \sum_{i=1}^N \partial_{x_{i}}\left(a_i(x,t,\nabla u) \right)  = f & {\rm in} \quad \Omega_{\infty} \equiv \Omega \times (0,+\infty) ,\\
u=0 & {\rm on} \quad \partial\Omega \times
(0,+\infty),\\
u(x,0)=u_0(x) & {\rm on} \quad \Omega,
\end{array}
\right.
\end{equation}
if  $u$ is  a  solution of \eqref{P_parab} for every $T>0$.
\end{definition}

We have the following results.

\begin{theorem}\label{thuni}
Let  $\Omega$ be a bounded domain of $\mathbb{R}^N$. Assume that \eqref{ii1}-\eqref{ii3} and \eqref{assupi} hold true and that it results   $u_0 \in L^1(\Omega)$ and  $f \in L^1(\Omega_T)$. Then there exists a weak solution $u$ of problem \eqref{P_parab}.
 Moreover, $u$ is in $C([0,T];L^1(\Omega))$ and  is the unique  solution of  \eqref{P_parab}  belonging to $C([0,T];L^1(\Omega))$ which is obtained by approximation (SOLA).

\noindent Finally, if $u$ is the unique solution of  \eqref{P_parab} obtained by approximation belonging to $C([0,T];L^1(\Omega))$, the following estimate holds true for every $ t \in (0,T)$
\begin{equation}\label{dipdati}
\|u(t)-v(t)\|_{ L^1(\Omega)} \leq \|u_0-v_0\|_{ L^1(\Omega)} + \|f-g\|_{ L^1(\Omega \times (0,t))}\,,
\end{equation}
where $v $ is the unique  solution obtained by approximation belonging to $C([0,T];L^1(\Omega))$ of the following problem
\begin{equation}\label{vanis}
 \left\{\begin{array}{ll}
v_t- \sum_{i=1}^N \partial_{x_{i}}\left(a_i(x,t, \nabla v ) \right)  = g & {\rm in} \quad \Omega_T ,\\
v=0 & {\rm on} \quad \partial\Omega \times
(0,T),\\
v(x,0)=v_0(x) & {\rm on} \quad \Omega,
\end{array}
\right.
\end{equation}
with $v_0 \in L^1(\Omega)$ and  $g \in L^1(\Omega_T)$.
\end{theorem}

\medskip
\begin{remark}
 We stress that the isotropic version of the previous result is proved in \cite{debmon} when $p_i \equiv p \quad \forall i=1,...,N$ .
\end{remark}

We point out that is also possible to prove the existence of a global solution of \eqref{P_parab}.
Indeed,  in the following results  we will show that there exists a global solution which is constructed by approximation, i.e. $u$ is a solution of \eqref{P_parab} obtained by approximation for every $T>0$ (according to Definition \ref{global}). In details we have the following result.
\begin{theorem}\label{glothuni}
 Let  $\Omega$ be a bounded domain of $\mathbb{R}^N$. Assume that \eqref{ii1}-\eqref{ii3} and \eqref{assupi}  hold true for every $T>0$.    If $u_0 \in L^1(\Omega)$ and  $f \in L^1_{loc}([0,\infty);L^1(\Omega))$, then there exists a global solution  $u$ of problem  \eqref{gloanis}.

\noindent Moreover, $u$ is in $C_{loc}([0,+\infty);L^1(\Omega))$ and  is the unique global  solution of   \eqref{gloanis} \nc obtained by approximation  belonging to $C_{loc}([0,+\infty);L^1(\Omega))$.

\noindent Finally, estimate \eqref{dipdati} holds true for every $t>0$, where now $v$ is the  unique global  solution of \eqref{vanis}  obtained by approximation  belonging to $C_{loc}([0,+\infty);L^1(\Omega))$.
\end{theorem}



Now let us replace  assumption \eqref{ii3} with the following stronger monotony assumption
\begin{eqnarray}
\label{ii4}&&     ( a_i(x,t, \xi )- a_i(x,t, \eta ) )\, (\xi_i - \eta_i ) \geq  \gamma |\xi_i - \eta_i|^{p_i}\,, \quad  \gamma > 0\,,
\end{eqnarray}
\noindent for a.e. $(x,t) \ \in \Omega_T$ and for any vector $\xi,\eta$  in $\mathbb R^N$.

We stress that in the case of the model problem \eqref{anis} the assumption \eqref{ii4} holds true if $p_{\min} = \min\{p_1,...,p_N\}\geq 2$.


As a first result we show that all the global  solutions obtained by approximation have the same asymptotic behavior, independently by the value of the initial datum.
In details, we have the following result.
\begin{theorem}\label{asinto}
Let  $\Omega$ be a bounded domain of $\mathbb{R}^N$.  Assume  $u_0 \in L^1(\Omega)$, $v_0 \in L^1(\Omega)$ and that \eqref{ii1},\eqref{ii2}, \eqref{assupi} and \eqref{ii4} hold true.

 If $u$ and $v$ are the unique  solutions in $C([0,T];L^1(\Omega))$ obtained by approximation of, respectively, \eqref{P_parab} and \eqref{vanis} with $g = f \in L^1(\Omega_T)$, then $u-v$ belongs to $L^{\infty}(\Omega \times (t,T)) \cap L^{\overrightarrow{p}}\left( t,T;W_{0}^{1,\overrightarrow
{p}}(\Omega)\right) $ and the following estimates hold true for every $t \in (0,T)$

\begin{equation}\label{uno}
\|u-v\|_{L^{\infty}(\Omega \times [t,T])} \leq C_1 \, \frac{\|u_0-v_0\|^{\frac{\overline{p}}{N(\overline{p}-2)+ \overline{p} }}_{L^{1}(\Omega)}}{t^{\frac{N}{N(\overline{ p}-2)+ \overline{p} }}} ,
\end{equation}
\begin{equation}\label{graduno}
\int_t^T\!\!\int_{\Omega}  \sum_{i=1}^N | \partial_{x_i}(u-v)|^{p_i} \leq  C_0 C_1^2 \, \frac{\|u_0-v_0\|^{\frac{2\overline{p}}{N(\overline{p}-2)+ \overline{p} }}_{L^{1}(\Omega)}}{t^{\frac{2N}{N(\overline{ p}-2)+ \overline{p} }}} ,
\end{equation}

\noindent where $C_0 = \frac{|\Omega|}{2\gamma}$ and $C_1 = C(N,\overline{p},\gamma)$. Moreover, it results
\begin{equation}\label{expo}
\|u-v\|_{L^{\infty}(\Omega \times [t,T])} \leq C_2 \, \frac{\|u_0-v_0\|_{L^{1}(\Omega)}}{t^{\frac{N}{ 2 }} e^{ \sigma t}},  \quad  \hbox{\rm if} \,\, \overline{ p} = 2\,,
\end{equation}
\begin{equation}\label{uni}
\|u-v\|_{L^{\infty}(\Omega \times [t,T])} \leq \frac{C_*}{t^{\frac{1}{\overline{p}-2}}} \,, \quad \hbox{\rm if} \,\,  2<\overline{ p}<N\,,
\end{equation}
and
\begin{equation}\label{gradexpo}
\int_t^T\!\!\int_{\Omega}  \sum_{i=1}^N | \partial_{x_i}(u-v)|^{p_i} \leq C_0 C_2^2 \, \frac{\|u_0-v_0\|^2_{L^{1}(\Omega)}}{t^N e^{ 2\sigma t}},  \quad  \hbox{\rm if} \,\, \overline{ p} = 2\,,
\end{equation}
\begin{equation}\label{graduni}
\int_t^T\!\!\int_{\Omega}  \sum_{i=1}^N | \partial_{x_i}(u-v)|^{p_i} \leq \frac{(C_*)^2C_0}{t^{\frac{2}{\overline{p}-2}}} \,, \quad \hbox{\rm if} \,\,  2<\overline{ p}<N\,,
\end{equation}
\noindent where $C_2=C_2(N,\gamma)$,  $\sigma$ is a positive constant depending only on $N$ and $|\Omega|$, and $C_*= C_*(N,\overline{p},|\Omega|)$. Finally, if $f$  belongs to $ L^1_{loc}([0,\infty);L^1(\Omega))$ and $u$ and $v$ are the unique global  solutions   obtained by approximation of, respectively, \eqref{gloanis} and \eqref{vanis} with again $g =f$, which belong to  $C_{loc}([0,+\infty);L^1(\Omega))$, then estimates \eqref{uno}-\eqref{graduni}  hold true for every $t >0$ and $T> t$. In particular it results
\begin{equation*}\label{lima}
\lim_{t \rightarrow +\infty}\|u-v\|_{L^{\infty}(\Omega \times [t,+\infty))} = 0\,.
\end{equation*}
\end{theorem}

\begin{remark}
Notice that  if $ 2<\overline{ p}<N$, combing estimates \eqref{uno} and \eqref{uni} it follows that
\begin{equation*}\label{zuni}
\|u-v\|_{L^{\infty}(\Omega \times [t,T])} \leq \min \left \{C_1 \frac{\|u_0-v_0\|^{\frac{\overline{p}}{N(\overline{p}-2)+ \overline{p} }}_{L^{1}(\Omega)}}{t^{\frac{N}{N(\overline{ p}-2)+ \overline{p} }}},\frac{C_*}{t^{\frac{1}{\overline{p}-2}}} \right\}\,,
\end{equation*}
for every $0<t<T$.
\end{remark}

\begin{remark}
 Theorem \ref{asinto} extends to anisotropic evolution problems what  is  known in the isotropic case (see \cite{fortmon}) and  it  is quite surprising because it states that the function $u-v$ belongs to  $L^{\infty}(\Omega \times (t,T)) \cap L^{\overrightarrow{p}}\left( t,T;W_{0}^{1,\overrightarrow
{p}}(\Omega)\right) $ even if both the solutions $u$ and $v$ are generally unbounded on the set $\Omega \times (t,T)$ and not belonging to the anisotropic space  $ L^{\overrightarrow{p}}\left( t,T;W_{0}^{1,\overrightarrow
{p}}(\Omega)\right) $.
\end{remark}

A consequence of the previous theorem is the following result.
\begin{theorem}\label{asiauto}
Let  $\Omega$ be a bounded domain of $\mathbb{R}^N$.  Assume  $u_0 \in L^1(\Omega)$ and that \eqref{ii1},\eqref{ii2},\eqref{assupi} and \eqref{ii4} hold true. We also require $f(x,t) \equiv f(x) \in L^1(\Omega)$ and  that it results
\begin{equation*}\label{iphauto}
a_i(x,t,\xi) = a_i(x,\xi) \quad \hbox{\rm for every} \quad i= 1, \cdots, N\,,
\end{equation*}
a.e. $(x,t) \in \Omega_T$ and for every $\xi \in \R^N$.
Then there exists a unique weak solution $w$  obtained by approximation of the following stationary problem  (see Definition \ref{defsolelli} below)
\begin{equation}\label{P_ell}
\left\{
\begin{array}
[c]{lll}%
\sum_{i=1}^N  \partial_{x_{i}}(a_i(x,\nabla w))=f(x) & \quad \hbox{\rm in }%
\Omega\,,  & \\
 w(x)=0 & \quad \hbox{\rm on }\;\partial\Omega\,. &
\end{array}
\right.
\end{equation}
Moreover, if $u$ is the unique global solution  in $C_{loc}([0,+\infty);L^1(\Omega))$ of \eqref{gloanis} which is obtained by approximation, then it results
\begin{equation}\label{limauto}
\lim_{t \rightarrow +\infty}\|u-w\|_{L^{\infty}(\Omega \times [t,+\infty))} = 0\,.
\end{equation}
\end{theorem}
We recall that by a weak solution of  \eqref{P_ell} obtained by approximation we mean the following.
\begin{definition}\label{defsolelli}
We say that $w \in W^{1,1}_0(\Omega)$ is a weak solution of  \eqref{P_ell} if
 $a_i(x,\nabla w) \in L^1(\Omega)$ for every $i=1,...,N$, and
\begin{equation*}
  \int_{\Omega} \sum_{i=1}^{N} a_i(x,\nabla w) \varphi_{ x_i}  \,dx =
  \int_{\Omega} f \varphi \,dx
\end{equation*}
for every $\varphi\in C^{\infty}_0(\overline{\Omega})$.

\noindent Moreover,
  a weak solution $w$ of  \eqref{P_ell} is a weak solution  of  \eqref{P_ell} obtained by approximation  if it is the a.e. limit in $\Omega$ of the solutions $w_n\in L^{\infty}(\Omega)\cap W_{0}%
^{1,\overrightarrow{p}}(\Omega)$ of the following (approximating) problems
\begin{equation}\label{elP_approx}
\left\{
\begin{array}
[c]{lll}%
\sum_{i=1}^N  \partial_{x_{i}}(a_i(x,\nabla w_n))=f_n(x) & \quad \hbox{\rm in }%
\Omega\,,  & \\
 w_n=0 & \quad \hbox{\rm on }\;\partial\Omega  , &
\end{array}
\right.  %
\end{equation}
where the data $f_n$  satisfy
\begin{equation}\label{apprell}
f_n\in L^{\infty}(\Omega), \quad f_n \rightarrow f \quad \text{in} \quad L^1(\Omega).
\end{equation}
\end{definition}

We show now that our equation has a very strong regularizing property.
In detail, we have the following result.

\begin{theorem}\label{regLinfty}
Let  $\Omega$ be a bounded domain of $\mathbb{R}^N$. Let $u_0$  be in $L^1(\Omega)$ and $f $ be in $L^m(\Omega_T)$, with $m > 1 + \frac{N}{\overline{p}}$. Assume that  \eqref{ii1},\eqref{ii2},\eqref{assupi} and \eqref{ii4} hold true.
 If   $u$ is the unique weak solution of \eqref{P_parab} in $C([0,T];L^1(\Omega))$ obtained by approximation then $u$ belongs to $L^{\infty}(\Omega \times (t,T))$, for every $t \in (0,T)$.

Moreover,  there exists a positive constant $C_3$ (see formula \eqref{defC3}) depending only on $\|f\|_{L^m(\Omega_T)}$, $N$, $p_i$, $t$, $\alpha$, $\|u_0\|_{L^{1}(\Omega)}$ and $|\Omega|$, such that the following estimate holds true
$$
\|u\|_{L^{\infty}(\Omega \times (t,T))} \leq  C_3\,.
$$
Furthermore, if $ \bar p > 2$, there exists a constant $C_4$  (see formula \eqref{defC4}) independent of $u_0$ (depending only on $\|f\|_{L^m(\Omega_T)}$, $N$, $p_i$, $t$, $\alpha$ and $|\Omega|$) such that the following universal bound is satisfied
$$
\|u\|_{L^{\infty}(\Omega \times (t,T))} \leq  C_4\,.
$$
Finally, it results
\begin{equation*}\label{reggrad}
u \in L^{\overrightarrow{p}}\left(  t,T;W_{0}^{1,\overrightarrow
{p}}(\Omega)\right) \,, \quad \quad \hbox{\rm for every} \quad t \in (0,T)\,.
\end{equation*}
\end{theorem}

\begin{remark}
We observe that if $u_0$ is only assumed to be in $L^1(\Omega)$, the   $L^{\infty}(\Omega \times (t,T))$-regularity  for every $t \in (0,T)$ of $u$  was already proved in  \cite{decay} in the particular case of a null forcing term $f$ in \eqref{P_parab}.
\end{remark}

More in general, we can prove the following  regularity property of the SOLA weak solutions $u$  which reveals that as soon as $t>0$ the regularity of the initial datum $u_0$ does not influence the regularity of $u$ since it earns, as soon as $t>0$, the same  regularity of a solution which has a bounded initial datum.
\begin{theorem}\label{regLrs}
Let  $\Omega$ be a bounded domain of $\mathbb{R}^N$. Let  $u_0$  be in $L^1(\Omega)$, $v_0 \in L^{\infty}(\Omega)$ and  $f \in L^m(\Omega_T)$, $m \geq 1$. Assume   that  \eqref{ii1},\eqref{ii2},\eqref{assupi} and \eqref{ii4} hold true.  If
  $u$  and $v$ are the unique weak solution in $C([0,T];L^1(\Omega))$ obtained by approximation of, respectively, \eqref{P_parab} and \eqref{vanis} with $g = f$,    then for every $t \in (0,T)$ the following property holds true	
	\begin{equation}\label{stessa}
	v \in L^r(t,T;L^s(\Omega)) \quad \Rightarrow \quad u \in L^r(t,T;L^s(\Omega))\,,
	\end{equation}
	where $r, s \in [1,+\infty]$.
	
	Moreover, if $v$ belongs to $L^{\overrightarrow{p}}\left(  t,T;W_{0}^{1,\overrightarrow
{p}}(\Omega)\right) $, with $t \in (0,T)$, then it results
\begin{equation}\label{regdu}
u \in L^{\overrightarrow{p}}\left( t,T;W_{0}^{1,\overrightarrow
{p}}(\Omega)\right)\,.
\end{equation}
\end{theorem}

\begin{remark}
We point out that the previous result extends to the anisotropic problem \eqref{P_parab} the known regularity result proved in \cite{fortmon} in the isotropic case $p = p_i$ for every $i=1,\cdots,N$ .
\end{remark}

\begin{remark}
We emphasize that a consequence of the previous results are some regularity results for the model problem \eqref{anis} when $u_0\in L^1(\Omega)$ and $p_{min} \geq 2$. More precisely, if $f \in L^m(\Omega_T)$ with $m > 1 + \frac{N}{\overline{p}}$ and $\overline{p} < N$, and if $u$ is the unique weak solution of \eqref{anis} in $C([0,T];L^1(\Omega))$ obtained by approximation then $u$ belongs to $L^{\infty}(\Omega \times (t,T))\cap L^{\overrightarrow{p}}\left(  t,T;W_{0}^{1,\overrightarrow
{p}}(\Omega)\right)$, for every $t \in (0,T)$.
Moreover, if $f \in L^m(\Omega_T)$, $m \geq 1$ and
  $u$  is the unique weak solution in $C([0,T];L^1(\Omega))$ obtained by approximation of \eqref{anis}    then
	\begin{equation*}\label{xstessa}
	v \in L^r(t,T;L^s(\Omega)) \quad \Rightarrow \quad u \in L^r(t,T;L^s(\Omega)) \quad \forall t \in (0,T)\,,
	\end{equation*}
	where $v$ is the unique weak solution in $C([0,T];L^1(\Omega))$ obtained by approximation of the following problem
	$$
 \left\{\begin{array}{ll}
v_t- \sum_{i=1}^N \partial_{x_i}\left( |v_{x_i}|^{p_i-2}v_{x_i} \right)  = f & {\rm in} \quad \Omega_T ,\\
v=0 & {\rm on} \quad \partial\Omega \times
(0,T),\\
v(x,0)=v_0(x) & {\rm on} \quad \Omega\,,
\end{array}
\right.
$$
with $v_0 \in L^{\infty}(\Omega)$.
\end{remark}

\section{Preliminaries }\label{secpre}

In this section, for the convenience of the reader, we recall some definition and results that we will use in the proofs.

\subsection{Anisotropic spaces. }\label{aniSpace}
Let $\Omega$ be an bounded open subset of $\mathbb{R}^{N}$, $N\geq2$, and let $1< p_{1},\ldots,p_{N}<\infty$ be
$N$ real numbers. We define the anisotropic Sobolev space $W_{0}^{1,p_{i}%
}(\Omega)$ as the closure of $C_{0}^{\infty}(\Omega)$ with respect to the
norm
\[
\left\Vert u\right\Vert _{W_{0}^{1,p_{i}}(\Omega)}=\left\Vert u\right\Vert
_{L^{1}(\Omega)}+\left\Vert u_{x_i}\right\Vert _{L^{p_{i}}(\Omega
)}.
\]
We recall that also in the anisotropic setting a Poincar\'{e}-type inequality holds
(see \cite{FGK}). More precisely if $u\in C_{0}^{\infty}(\Omega)$ with $\Omega$  a open bounded set with Lipschitz
continuous boundary there exists a constant $C_P$, depending only on the diameter of $\Omega$,
such that
\begin{equation}\label{dis poincare}
\left\Vert u\right\Vert _{L^{p_i}(\Omega)}\leq\ C_P\;\left\Vert u_{x_i}\right\Vert _{L^{q}(\Omega)} \quad \quad i=1,...,N.%
\end{equation}
We set $W_{0}^{1,\overrightarrow{p}}(\Omega)=\displaystyle\bigcap_{i=1}%
^{N}W_{0}^{1,p_{i}}(\Omega)$ with the norm
\begin{equation*}
\left\Vert u\right\Vert _{W_{0}^{1,\overrightarrow{p}}(\Omega)}=\overset
{N}{\underset{i=1}{\sum}}\left\Vert u_{x_i}\right\Vert _{L^{p_{i}%
}(\Omega)}\,, \label{Sob_norm}%
\end{equation*}
and we denote its dual by $\left(  W_{0}^{1,\overrightarrow{p}}(\Omega
)\right)  ^{\prime}.$
Here, for every $t \in [0,T)$, we define
$$
L^{\overrightarrow{p}}\left(  t,T;W_{0}^{1,\overrightarrow
{p}}(\Omega)\right)  =\displaystyle\bigcap_{i=1}^{N}L^{p_{i}}\left(
t,T;W_{0}^{1,p_{i}}(\Omega)\right)\,,
$$
where
$$
L^{p_i}(t,T;W_0^{1,p_i}(\Omega)) = \{ v \in L^{1}(t,T;W_0^{1,1}(\Omega)): v_{x_i}\in L^{p_i}(\Omega \times (t,T))\}\,.
$$


We recall that if  $u\in C_0^{1}(\mathbb{R^N})$ it is well known that if  $1 < \bar{p}<N$ the following anisotropic Sobolev inequality holds true (see Theorem 1.2 in \cite{tro})
\begin{equation}\label{Sob}
\|u\|_{L^{\bar{p}^*}(\mathbb{R^N})} \le C_S\prod_{i=1}^N \|  u_{x_i}\|_{L^{p_i}(\mathbb{R^N})}^{\frac 1N},
\end{equation}
where $C_S$ is a constant depending only on $N$,  $\bar p$ as in \eqref{p bar} and $\bar{p}^*=\frac{N\bar{p}}{N-\bar p}\,$. Using the inequality between geometric and arithmetic mean we can replace the right-hand-side of \eqref{Sob} with $ \sum_{i=1}^N \|  u_{x_i}\|_{L^{p_i} }$.
%
When $\overline{p}<N$ and $\Omega$ is a bounded open set with Lipschitz boundary, the Sobolev inequality
implies the continuous embedding of the space $W_{0}^{1,\overrightarrow{p}%
}(\Omega)$ into $L^{q}(\Omega)$ for every $q\in\lbrack1,\bar{p}^{\ast}]$. On
the other hand, the Poincar\'{e} inequality \eqref{dis poincare} assure the continuity of the embedding $W_{0}^{1,\overrightarrow{p}%
}(\Omega)\subset L^{p_{\max}}(\Omega)$ with $p_{\max}:=\max\{p_{1},\ldots,p_{N}\}$. If the exponents $p_{i}$ are not closed enough it may happen that $\bar{p}^{\ast}<p_{\max}$ , then in this case $p_{\infty}:=\max\{\bar{p}^{\ast},p_{\max}\}$ turns out to be the critical exponent in the anisotropic Sobolev embedding.

\noindent

\bigskip
\subsection{Decay estimates}

In this subsection we recall some
 results about decay estimates which can be deduced simply by integral estimates.

\begin{theorem}[Theorem 2.1 in \cite{decay}]\label{teo3}
Assume that
\begin{equation*}\label{somu}
w \in C((0,T);L^{r}(\Omega)) \cap L^{b}(0,T;L^q(\Omega)) \cap C([0,T);L^{r_0}(\Omega))
\end{equation*}
 where $\Omega$ is an open  set of $\R^N$, (not necessary bounded), $N \geq 1$, $0<T\leq +\infty$ and
\begin{equation*}\label{rqte}
1 \leq r_0 < r <q \leq + \infty, \quad b_0 < b < q, \quad b_0 = \frac{(r-r_0)}{1-\frac{r_0}{q}}\,.
\end{equation*}
Suppose that $w$   satisfies the following integral estimates  for every \, $k>0$,
\begin{equation}\label{esqte}
\int_\Omega |G_k(w)|^{r}(t_2) dx - \int_\Omega |G_k(w)|^{r}(t_1) dx + c_1 \int_{t_1}^{t_2} \| G_k(w)(\tau)\|^b_{L^q(\Omega)}  d\tau \leq 0
\end{equation}
  for every  $0 < t_1 < t_2 < T$  and
\begin{equation}\label{esr0}
  \|G_k(w)(t)\|_{L^{r_0}(\Omega)} \leq c_2 \|G_k(w)(t_0)\|_{L^{r_0}(\Omega)} \quad \hbox{\rm for every} \,\,\, 0 \leq t_0 < t < T,
\end{equation}
where $c_1$ and $c_2$  are positive constants independent of $k$. Finally, let us define
\begin{equation}\label{ipu0}
   w_0 \equiv w(x,0) \in L^{r_0}(\Omega).
\end{equation}
Then there exists a positive constant $C_1$ (see formula (4.11) in \cite{decay}) depending only on $N$, $c_1$, $c_2$, $r$, $r_0$, $q$ and $b$ such that
\begin{equation*}\label{tesi}
\| w(t)\|_{L^\infty(\Omega)} \leq C_1 \frac{\|w_0\|^{h_0}_{L^{r_0}(\Omega)}}{t^{h_1}} \quad  \hbox{\rm for every} \,\,\, t \in (0,T),
\end{equation*}
where
\begin{equation}\label{deh01}
h_1 = \frac{1}{b - (r- r_0)- \frac{r_0b}{q}}, \quad h_0 = h_1 \left(1-\frac{b}{q}\right)r_0.
\end{equation}
\end{theorem}
If  $\Omega$ has finite measure we have an exponential decay if  $b = r$ and universal bounds if  $b > r$. More in detail
we have the following  result.

\begin{theorem}[Theorem 2.2 in \cite{decay}]\label{genelim} Let the assumptions of Theorem \ref{teo3} hold true.

\noindent If\ $\Omega$ has finite measure and $b = r$
 the following exponential decay occurs
\begin{equation*}\label{tesexp}
\| w(t)\|_{L^\infty(\Omega)} \leq C_2 \frac{\|w_0\|_{L^{r_0}(\Omega)}}{t^{h_1}e^{\sigma t}} \quad  \hbox{\rm for every} \,\,\, t \in (0,T),
\end{equation*}
where  $C_2$  is a positive constant depending only on $N$, $c_1$, $c_2$, $r$, $r_0$ and $q$ (see formula (4.17) in \cite{decay}), $u_0$ is as in \eqref{ipu0}, $h_1$ is as in \eqref{deh01} with $r=b$, i.e.
\begin{equation*}\label{defHi}
h_1 =  \frac{1}{ r_0 \left( 1 - \frac{r}{q}\right)},
\end{equation*}
and
\begin{equation*}\label{defsig}
\sigma = \frac{c_1 \kappa}{4(r-r_0)|\Omega|^{1-\frac{r}{q}}}, \quad \kappa \quad \hbox{\rm arbitrarily fixed in} \quad \left(0,1-\frac{r_0}{ r}\right),
\end{equation*}
where $|\Omega|$ denotes the measure of $\Omega$.
If otherwise  $\Omega$ has finite measure and $ b > r$ we have the following universal bound
\begin{equation*}\label{teuniv}
\| w(t)\|_{L^\infty(\Omega)} \leq  \frac{C_{*}}{t^{h_2}} \quad \hbox{\rm for every} \,\,\, t \in (0,T),
\end{equation*}
where
\begin{equation*}\label{defhu}
h_2= h_1 + \frac{h_0}{b -r} = \frac{1}{b-r},
\end{equation*}
 and $C_{*} $ is a constant depending only on $r$, $r_0$, $q$, $b$, $c_1$, $c_2$ and the measure of $\Omega$ (see formula (4.19) in \cite{decay}).
\end{theorem}

\section{Proof of the results}\label{secproofs} In this section we prove all the results stated in Section \ref{secstate}.

\subsection{Proof of Theorem \ref{thuni}}

The existence of a solution constructed by approximation for the anisotropic parabolic problem \eqref{P_parab} according to the above Definition \ref{defsol} is proved in \cite{Li}.

We observe that thanks to assumption \eqref{ii3} we can adapt the proof in \cite{lower} obtaining  that every approximating sequence $u_n$ of \eqref{P_approx} is a Cauchy sequences in $C([0,T];L^1(\Omega))$ and hence every solution $u$ obtained by approximation belongs to $ C([0,T];L^1(\Omega))$ and is the limit  in $C([0,T];L^1(\Omega))$ of the sequence $u_n$. In detail, let $\delta > 1$ and  $\varphi =  \left\{ 1 - \frac{1}{[1+|u_n-u_m|]^{\delta}}\right\}$ sign$(u_n-u_m)$. Choosing  $\varphi$ as test function in the equation satisfied by $u_n$ and in that satisfied by $u_m$ and subtracting the results  we get (for every $t \in [0,T]$)

\begin{eqnarray}
&& \int_{\Omega}|u_n(t)-u_m(t)|+\frac{1}{1-\delta}\int_{\Omega}\left(1-\frac{1}{(1+|u_n(t)-u_m(t)|)^{\delta-1}}\right) \nonumber\\
&& +\delta \int_0^t\int_{\Omega} \sum_{i=1}^{N}(a_i(x,t, \nabla u_n)-a_i(x,t, \nabla u_m))\frac{\partial_{x_i}(u_n-u_m)}{(1+|u_n(t)-u_m(t)|)^{\delta+1}}\leq \\
&& +\int_{\Omega}|u_{0,n}-u_{0,m}|+\frac{1}{1-\delta}\int_{\Omega}\left(1-\frac{1}{(1+|u_{0,n}-u_{0,m}|)^{\delta-1}}\right)+\int_0^t\int_{\Omega}|f_n-f_m|. \nonumber
\end{eqnarray}
By \eqref{ii3} and the previous inequality  we obtain
\begin{align*}
\int_{\Omega}|u_n(t)-u_m(t)|\leq \frac{|\Omega|}{\delta-1}+\int_{\Omega}|u_{0,n}-u_{0,m}|+\int_0^t\int_{\Omega}|f_n-f_m|.
\end{align*}
Finally letting $\delta \rightarrow +\infty$ we deduce
\begin{equation*}\label{didi}
\int_{\Omega}|u_n(t)-u_m(t)|\leq \int_{\Omega}|u_{0,n}-u_{0,m}|+\int_0^t\int_{\Omega}|f_n-f_m|\,.
\end{equation*}
from which it follows that $u_n$ is a Cauchy sequences in $C([0,T];L^1(\Omega))$ being  $u_{0,n}$ and $f_n$  Cauchy sequences in, respectively,  $L^1(\Omega)$ and  $L^1(\Omega_T)$.

We show now that if $u$ and $v$ belong to   $C([0,T];L^1(\Omega))$ and are  solutions constructed by approximation of, respectively,  \eqref{P_parab} and \eqref{vanis}, then estimate \eqref{dipdati} holds true and consequently it is unique the solution  $u$  of \eqref{P_parab} belonging to $C([0,T];L^1(\Omega))$ which is constructed by approximation.

In order to prove \eqref{dipdati},
let $u$ be the a.e. limit in $\Omega_T$ of the solution $u_n$ of the approximating problems \eqref{P_approx}  and $v$ be the a.e. limit in $\Omega_T$ of the solution $v_n\in L^{\infty}(\Omega_T)\cap C([0,T];L^2(\Omega))\cap L^{\overrightarrow{p}}\left(  0,T;W_{0}%
^{1,\overrightarrow{p}}(\Omega)\right)$ of the following approximating problems
\begin{equation}\label{P_approxv}
\left\{
\begin{array}
[c]{lll}%
(v_n)_{t}+\sum_{i=1}^N  \partial_{x_i}(a_i(x,t,\nabla v_n))=g_n(x,t) & \quad \hbox{\rm in }%
\Omega_{T},  & \\
 v_n=0 & \quad \hbox{\rm on }\;\partial\Omega\times\left(  0,T\right)  , &\\
v_n(x,0)=v_{0,n}(x) & \quad \hbox{\rm on }\;\Omega, &
\end{array}
\right.  %
\end{equation}
where $g_n$ and $v_{0,n}$ satisfy
\begin{equation*}
g_n\in L^{\infty}(\Omega_T), \quad g_n\rightarrow g \quad \text{in} \quad L^1(\Omega_T),
\end{equation*}
\begin{equation*}
v_{0,n}\in L^{\infty}(\Omega), \quad v_{0,n}\rightarrow v_0 \quad \text{in} \quad L^1(\Omega).
\end{equation*}
 For what proved above we now that both the sequences $u_n$ and $v_n$ are Cauchy sequences in $C([0,T];L^1(\Omega))$ converging to, respectively, $u$ and $v$ in $C([0,T];L^1(\Omega))$.

Reasoning as before one can choose  $\varphi=\left(1-\frac{1}{(1+|u_n-v_n|)^{\delta}}\right) \rm{sign} $ $(u_n-v_n)$ as test function in  \eqref{P_approx} and in \eqref{P_approxv} with $\delta>1$ arbitrarily fixed, obtaining
%
%
\begin{align*}
\int_{\Omega}|u_n(t)-v_n(t)|\leq \frac{|\Omega|}{\delta-1}+\int_{\Omega}|u_{0,n}-v_{0,n}|+\int_0^t\int_{\Omega}|f_n-g_n|.
\end{align*}

Finally letting $\delta \rightarrow +\infty$, since $u_n$ and $v_n$ are Cauchy sequence in $C([0,T];L^1(\Omega))$ it follow that
\begin{equation*}
\int_{\Omega}|u(t)-v(t)|\leq \int_{\Omega}|u_{0}-v_{0}|+\int_0^t\int_{\Omega}|f-g|.
\end{equation*}

\qed

\subsection{Proof of Theorem \ref{glothuni}}

By Theorem \ref{thuni} there exists a unique solution of \eqref{P_parab} obtained by approximation which belongs to $C([0,T];L^1(\Omega))$. For any $T>0$ fixed, we define $u(x,t)$ in $\Omega_T$ such a solution. Notice that the definition is well posed thanks to the uniqueness result in  Theorem \ref{thuni}.  By construction $u$ is defined in $\Omega \times (0,\infty)$  and  is  a solution obtained by approximation for every arbitrary fixed $T>0$  which   belongs to  $C([0,T];L^1(\Omega))$ for every $T>0$. Consequently $u$ is in $C_{loc} ([0,\infty);L^1(\Omega))$ and is a global solution obtained by approximation. Moreover, $u$ is also unique (again by Theorem \eqref{thuni}). Finally, by construction and thanks to Theorem \ref{thuni} estimate \eqref{dipdati} holds for every $t>0$.

\qed

\subsection{Proof of Theorem \ref{asinto}}


Let $u$ and $v$ be the unique solution in $C([0,T];L^1(\Omega))$ obtained by approximation of, respectively, \eqref{P_parab} and \eqref{vanis} with $g = f$, whose existence and uniqueness is guaranteed by Theorem \ref{thuni}. Using assumption \eqref{ii3}, which is implied by assumption \eqref{ii4}, we get that $u$ and $v$ are limits of approximating sequences $u_n$ of \eqref{P_approx} and $v_n$ of \eqref{P_approxv} with $g_n = f_n$, respectively, which are Cauchy sequences in $C([0,T];L^1(\Omega))$.  We recall that being $g = f$ and since the solutions in $C([0,T];L^1(\Omega))$ of \eqref{P_parab} and \eqref{vanis} obtained by approximation are unique, we can choose in the approximating problems  \eqref{P_approx} and  \eqref{P_approxv} the same approximation of the data  $g_n = f_n$.

Let $k >0$ and choose as test function $\varphi = G_k(u_n-v_n) \equiv \left(|u_n-v_n|-k\right)_+$sign$(u_n-v_n)$ in \eqref{P_approx} and in \eqref{P_approxv}  with $g_n = f_n$ the same approximation chosen in \eqref{P_approx} (recall that here  $g = f$). Subtracting the equations obtained in this
 way and using \eqref{ii4}  we get for every $0 < t_1 < t_2 \leq T$
\begin{equation}\label{stima G}
  \frac{1}{2}\int_{\Omega}|G_k(u_n-v_n)|^2(t_2)-\frac{1}{2}\int_{\Omega}|G_k(u_n-v_n)|^2(t_1) + \gamma \int_{t_1}^{t_2}\!\!\int_{\Omega} \sum_{i=1}^{N}|\de_{x_i}G_k(u_n-v_n)|^{p_i}\leq 0.
\end{equation}
We recall that being $\bar p < N$ we can use \eqref{Sob} obtaining
$$
  \|G_k(u_n-v_n)\|_{L^{{\bar p}^*}(\Omega)} \leq C_S \prod_{j=1}^ N \left(\int_{\Omega} |\de_{x_j}G_k(u_n-v_n)|^{p_j}\right)^{\frac{1}{Np_j}}\leq C_S \left(\int_{\Omega}\sum_{i=1}^N|\de_{x_i}G_k(u_n-v_n)|^{p_i}\right)^{\frac{1}{\bar p}}
$$
that is
\begin{equation}\label{G}
\|G_k(u_n-v_n)\|^{\bar p}_{L^{{\bar p}^*}(\Omega)}  \leq C_S^{\bar p} \int_{\Omega}\sum_{i=1}^N|\de_{x_i}G_k(u_n-v_n)|^{p_i}\,,
\end{equation}
where $C_S$ is a positive constant depending only on $N$ an $\bar p$.
By \eqref{G} and \eqref{stima G} we deduce
\begin{equation*}\label{2.3}
  \int_{\Omega}|G_k(u_n-v_n)|^2(t_2)-\int_{\Omega}|G_k(u_n-v_n)|^2(t_1)+\frac{2 \gamma} {C_S^{\bar p}}\int_{t_1}^{t_2}\|G_k(u_n-v_n)\|^{\bar p}_{L^{{\bar p}^*}(\Omega)}\leq 0{\color{red},}{\color{blue}.}
\end{equation*}

Hence the integral estimate \eqref{esqte} of Theorem \ref{teo3} is satisfied  by $w = u_n-v_n$ with $r = 2$, $c_1 = \frac{2 \gamma} {C_S^{\bar p}}$, $q={\bar p}^* $ and $b = \bar p $.
Now we want to prove also the other integral estimate \eqref{esr0} of Theorem \ref{teo3} is satisfied. To this aim we
 reason as in Theorem \ref{thuni}  choosing as test function in the approximating problems \eqref{P_approx} and \eqref{P_approxv} (with $g_n = f_n$) the function $\varphi=\left(1-\frac{1}{(1+|G_k(u_n-v_n)|)^{\delta}}\right)$sign$(u_n-v_n)$, for a suitable $\delta>1$. Subtracting the equations obtained in this way we get for every $0 < t_0 < t < T$
\begin{align}\label{Abis}
\int_{\Omega}&|G_k(u_n-v_n)|(t)+\frac{1}{1-\delta}\int_{\Omega}\left(1-\frac{1}{(1+|G_k (u_n-v_n)|(t))^{\delta-1}}\right)\\
\nonumber &+\delta \int_{t_0}^t\int_{|u_n-v_n|>k} \sum_{i=1}^{N}(a_i(x,t, \nabla u_n)-a_i(x,t, \nabla v_n)\frac{\de_{x_i}(u_n-v_n)}{(1+|G_k(u_n-v_n)|(t))^{\delta+1}}\leq \\
\nonumber &+\int_{\Omega}|G_k(u_n-v_n)|(t_0)+\frac{1}{1-\delta}\int_{\Omega}\left(1-\frac{1}{(1+|G_k(u_n-v_n)|(t_0))^{\delta-1}}\right).
\end{align}
By \eqref{ii3} (which is implied by assumption \eqref{ii4}) and inequality \eqref{Abis} we deduce
\begin{align*}
\int_{\Omega}|G_k(u_n-v_n)|(t)\leq \frac{|\Omega|}{\delta-1}+\int_{\Omega}|G_k(u_n-v_n)|(t_0).
\end{align*}
Thus, letting $\delta$ to infinite, 
 we obtain
\begin{align*}
\int_{\Omega}|G_k(u_n-v_n)|(t)\leq \int_{\Omega}|G_k(u_n-v_n)|(t_0)\,,
\end{align*}
and hence \eqref{esr0} is satisfied with $r_0 = c_2 = 1$.


 As observed in Remark \ref{rempbar}, the  assumption \eqref{assupi} implies  $\bar p> \frac{2N}{N+1}$, so  all the assumptions of Theorem \ref{teo3} are satisfied and hence the following estimate holds true
\begin{equation}\label{primadis2.7}
\|u_n(t)-v_n(t)\|_{L^{\infty}(\Omega)} \leq C_1 \frac{\|u_{0,n}-v_{0,n}\|^{\frac{\overline{p}}{N(\overline{p}-2)+ \overline{p} }}_{L^{1}(\Omega)}}{t^{\frac{N}{N(\overline{p}-2)+ \overline{p} }}}\,,
\end{equation}
where $C_1$ is a positive constant depending only on  $N$,  $\bar p$ and $\gamma$.
Notice that if  $\bar p > 2$, that means  $b > r$ in \eqref{esqte}, we can also apply Theorem \ref{genelim} obtaining in this particular case the following estimate
\begin{equation}\label{seconda2.7}
\|u_n(t)-v_n(t)\|_{L^{\infty}(\Omega)} \leq \frac{C_*}{t^{\frac{1}{\overline{p}-2}}},
\end{equation}
where $C_*$ is a positive constant depending only on $N$, $\overline{p}$ and $|\Omega|$.

Moreover, if  $\bar p =2$ it results $b = r$ in \eqref{esqte}  and hence  applying Theorem \ref{genelim} we get the following estimate
\begin{equation}\label{expon}
\|u_n(t)-v_n(t)\|_{L^{\infty}(\Omega)}  \leq C_2 \, \frac{\|u_{0,n}-v_{0,n}\|_{L^{1}(\Omega)}}{t^{\frac{N}{ 2 }} e^{ \sigma t}},
\end{equation}
 with $C_2$  a positive constant depending only on $N$ and $\gamma$ and $\sigma$ a positive constant depending only on $N$ and $|\Omega|$.

Now, recalling that $u-v$ is the a.e. limit in $\Omega_T$ of the sequence $u_n-v_n$ and that by construction $u_{0,n}-v_{0,n}$ converges to $u_0-v_0$ in $L^{1}(\Omega)$,  by the above $L^{\infty}$-estimates  \eqref{primadis2.7}-\eqref{expon} it follows that  the bound \eqref{uno}, \eqref{expo} and \eqref{uni} hold true. Consequently,  $u-v$ belongs to $L^{\infty}(\Omega \times (t,T))$ for every $t \in (0,T)$.

To conclude the proof, we show now that also the estimates \eqref{graduno}, \eqref{gradexpo} and \eqref{graduni} are satisfied and hence $u-v$ belongs also to $L^{\overrightarrow{p}}\left(  t,T;W_{0}^{1,\overrightarrow
{p}}(\Omega)\right) $, for every $t \in (0,T)$. To this aim we notice that by \eqref{stima G} (applied with  $0 < t_1 = t < t_2 = T$) we deduce that
\begin{equation}\label{zstima G}
    \int_{t}^{T}\!\!\int_{\Omega} \sum_{i=1}^{N}|\de_{x_i}G_k(u_n-v_n)|^{p_i}\leq \frac{1}{2  \gamma}\int_{\Omega}|G_k(u_n-v_n)|^2(t) \leq C_0 \|u_n(t)-v_n(t)\|^2_{L^{\infty}(\Omega)}\,,
\end{equation}
where $C_0 = \frac{|\Omega|}{2  \gamma}$. Thus, by \eqref{zstima G} and   \eqref{primadis2.7} we deduce
\begin{equation}\label{zgraduno}
\int_t^T\!\!\int_{\Omega}  \sum_{i=1}^N | \partial_{x_i}(u_n-v_n)|^{p_i} \leq  C_0C_1^2 \, \frac{\|u_{0,n}-v_{0,n}\|^{\frac{2\overline{p}}{N(\overline{p}-2)+ \overline{p} }}_{L^{1}(\Omega)}}{t^{\frac{2N}{N(\overline{ p}-2)+ \overline{p} }}} \,.
\end{equation}
Moreover, in the particular cases that $\overline{ p} \geq 2$, by  \eqref{zstima G}  and \eqref{seconda2.7},\eqref{expon} it follows
\begin{equation*}\label{zgradexpo}
\int_t^T\!\!\int_{\Omega}  \sum_{i=1}^N | \partial_{x_i}(u_n-v_n)|^{p_i} \leq C_0C_2^2 \, \frac{\|u_{0,n}-v_{0,n}\|^2_{L^{1}(\Omega)}}{t^N e^{ 2\sigma t}},  \quad  \hbox{\rm if} \,\, \overline{ p} = 2\,,
\end{equation*}
\begin{equation}\label{zgraduni}
\int_t^T\!\!\int_{\Omega}  \sum_{i=1}^N | \partial_{x_i}(u_n-v_n)|^{p_i} \leq \frac{(C_*)^2C_0}{t^{\frac{2}{\overline{p}-2}}} \,, \quad \hbox{\rm if} \,\,  2<\overline{ p}<N\,.
\end{equation}
Thus, by \eqref{zgraduno},\eqref{zgraduni} we can conclude that estimates \eqref{graduno}, \eqref{gradexpo} and \eqref{graduni} hold true.


\qed

\subsection{Proof of Theorem \ref{asiauto}}

 By  \cite{castro}  (see the proofs of Theorems 2.3 and 2.4 and Remark 3.7 in \cite{castro})  there exists at least  a weak solution  $ w(x)$ of \eqref{P_ell} obtained by approximation.  Hence,  $w$ is the a.e. limit in $\Omega$ of the solutions $w_n\in L^{\infty}(\Omega)\cap W_{0}%
^{1,\overrightarrow{p}}(\Omega)$ of \eqref{elP_approx} with the data $f_n$  satisfying
\eqref{apprell}.  Moreover  following the proofs of Theorems 2.3 and 2.4 in \cite{castro}, one get that
\begin{equation*}\label{convwnL1}
w_n \,\, \longrightarrow \,\, w \quad \hbox{\rm in } \,\, L^1(\Omega)\,.
\end{equation*}

We observe that  the following functions
 $$
w_n(x,t) \equiv  w_n(x)  \quad  \hbox{\rm for every} \,\, t >0\,{\color{blue},}
$$
belong to  $C_{loc}([0,+\infty];L^1(\Omega))$ and are global weak solutions of the following parabolic problems
$$
\left\{
\begin{array}
[c]{lll}%
(w_n)_{t}+\sum_{i=1}^N   \partial_{x_{i}}(a_i(x,\nabla w_n))=f_n(x) & \quad\text{in }%
\Omega_{T},  & \\
 w_n(x,t)=0 & \quad\text{on}\;\partial\Omega\times\left(  0,T\right)  , &\\
w_n(x,0)=w_n(x) & \quad\text{in}\;\Omega\,. &
\end{array}
\right.  %
$$
Hence,
the function
$$
w(x,t) \equiv w(x) \quad \hbox{\rm for every} \,\, t >0\,,
$$  belongs to $C_{loc}([0,+\infty];L^1(\Omega))$ and  is a global  solution  obtained by approximation of the following parabolic problem
\begin{equation*}\label{P_parabw}
\left\{
\begin{array}
[c]{lll}%
w_{t}+\sum_{i=1}^N   \partial_{x_{i}}(a_i(x,\nabla w))=f(x) & \quad\text{in }%
\Omega_{T}\,,  & \\
 w(x,t)=0 & \quad\text{on}\;\partial\Omega\times\left(  0,T\right)  , &\\
w(x,0)=w(x) & \quad\text{in}\;\Omega\,. &
\end{array}
\right.  %
\end{equation*}
Thus, by  Theorem \ref{glothuni} it follows that $w$ is the unique global solution in  $C_{loc}([0,+\infty];L^1(\Omega))$  obtained by approximation of problem   \eqref{vanis} with data $g(x,t) \equiv f(x)$ and $v_0(x) \equiv w(x)$.  Consequently,
 estimate \eqref{limauto} follows by Theorem \ref{asinto}.

To conclude the proof, it remains to show that $w$ is also the unique solution obtained by approximation of the stationary problem \eqref{P_ell} (according to Definition \ref{defsolelli}).

If $W$ is another solution obtained by approximation of the stationary problem \eqref{P_ell}, by what proved above we know that the following estimates hold true
\begin{equation}\label{wlimauto}
\lim_{t \rightarrow +\infty}\|u-w\|_{L^{\infty}(\Omega \times [t,+\infty))} = 0\,,
\end{equation}
\begin{equation}\label{WWlimauto}
\lim_{t \rightarrow +\infty}\|u-W\|_{L^{\infty}(\Omega \times [t,+\infty))} = 0\,.
\end{equation}
Notice that for every $t>0$ it results for a.e. $x \in \Omega$
$$
|w(x)-W(x)| \leq |w(x)-u(x,t)| + |u(x,t)-W(x)|  \leq \|u-w\|_{L^{\infty}(\Omega \times [t,+\infty))} + \|u-W\|_{L^{\infty}(\Omega \times [t,+\infty))}
$$
and hence letting $t \rightarrow +\infty$ in the previous inequality, by \eqref{wlimauto} and \eqref{WWlimauto} it follows that $w(x)=W(x)$ a.e. $x \in \Omega$.

\qed

\subsection{Proof of Theorem \ref{regLinfty}}

Let $v \in C([0,T];L^1(\Omega))$ be the unique solution  obtained by approximation of the following problem
\begin{equation*}\label{vP_parab}
\left\{
\begin{array}
[c]{lll}%
v_{t}+\sum_{i=1}^N \partial_{x_{i}}(a_i(x,t,\nabla v))=f(x,t) & \quad\text{in }%
\Omega_{T},  & \\
 v(x,t)=0 & \quad\text{on}\;\partial\Omega\times\left(  0,T\right)  , &\\
v(x,0)= 0 & \quad\text{in}\;\Omega\,. &
\end{array}
\right.  %
\end{equation*}
Since $v(x,0)= 0 \in L^{\infty}(\Omega)$ and by assumption $f \in L^m(\Omega_T)$ with  $m > 1 + \frac{N}{\overline{p}}$
we can apply Theorem 2.1 in \cite{P3} obtaining that $v \in L^{\infty}(\Omega_T)$ and satisfies
\begin{equation*}\label{stimainfv}
\|v\|_{L^{\infty}(\Omega_T)} \leq c,
\end{equation*}
where $c$ is a positive constant depending only on $\|f\|_{L^m(\Omega_T)}$, $N$, $p_i$ and $\alpha$. Thus, by Theorem \ref{asinto} we deduce that $u \in L^{\infty}(\Omega \times (t,T))$ for every $t \in (0,T)$. Indeed, thanks to the bound \eqref{uno} we deduce
\begin{equation}\label{defC3}
\|u\|_{L^{\infty}(\Omega \times (t,T))} \leq \|u-v\|_{L^{\infty}(\Omega\times (t,T))} + \|v\|_{L^{\infty}(\Omega\times (t,T))} \leq C_3 \equiv C_1 \, \frac{\|u_0\|^{\frac{\overline{p}}{N(\overline{p}-2)+ \overline{p} }}_{L^{1}(\Omega)}}{t^{\frac{N}{N(\overline{ p}-2)+ \overline{p} }}} + c\,.
\end{equation}
Notice that the constant $C_3$ depends only on  $\|f\|_{L^m(\Omega_T)}$, $N$, $p_i$, $t$, $\alpha$, $\|u_0\|_{L^{1}(\Omega)}$ and $|\Omega|$.

Furthermore, if $\bar p > 2$, proceeding exactly as above but replacing \eqref{uno} with \eqref{uni} we obtain the following universal estimate
\begin{equation}\label{defC4}
\|u\|_{L^{\infty}(\Omega \times (t,T))} \leq C_4 \equiv \frac{C_*}{t^{\frac{1}{\overline{p}-2}}} + c\,,
\end{equation}
where  $C_*=C_*(N,\bar p, |\Omega|)$  and consequently $C_4$ is independent of $u_0$ and depends only on  $\|f\|_{L^m(\Omega_T)}$, $N$, $p_i$, $t$, $\alpha$ and $|\Omega|$.

Finally,  running the proof of Theorem 2.1 in \cite{P3} one get that $v$ belongs to $ L^{\overrightarrow{p}}\left( 0,T;W_{0}^{1,\overrightarrow
{p}}(\Omega)\right)$, on the other hand  by Theorem \ref{asinto} we know that $u-v$ is in $L^{\overrightarrow{p}}\left( t,T;W_{0}^{1,\overrightarrow
{p}}(\Omega)\right)$, for every $t \in (0,T)$,  and so  we can conclude that  $u \in L^{\overrightarrow{p}}\left(  t,T;W_{0}^{1,\overrightarrow
{p}}(\Omega)\right)$.

\qed

\subsection{Proof of Theorem \ref{regLrs}}

Let $v \in C([0,T];L^1(\Omega))$ be the unique solution   obtained by approximation of \eqref{vanis} with $g \equiv f$.
The regularity result \eqref{stessa} is an immediate consequence of Theorem \ref{asinto}. As a matter of fact by \eqref{uno} we deduce that  the following bound is satisfied
\begin{equation}\label{star}
|u(x,t)| \leq \|u(t)-v(t)\|_{L^{\infty}(\Omega)} + |v(x,t)| \leq C_1 \, \frac{\|u_0-v_0\|^{\frac{\overline{p}}{N(\overline{p}-2)+ \overline{p} }}_{L^{1}(\Omega)}}{t^{\frac{N}{N(\overline{ p}-2)+ \overline{p} }}}  + |v(x,t)|\,,
\end{equation}
and hence if $v$ belongs to $L^r(t,T;L^s(\Omega))$, $r, s \in [1,+\infty]$, also  $ u $ is in $ L^r(t,T;L^s(\Omega))$.
Notice that if  $\bar p > 2$, proceeding as in \eqref{star} but replacing \eqref{uno} with \eqref{uni} we deduce that
$$
|u(x,t)| \leq \|u(t)-v(t)\|_{L^{\infty}(\Omega)} + |v(x,t)| \leq \frac{C_*}{t^{\frac{1}{\overline{p}-2}}} + |v(x,t)|\,.
$$
We point out that
the previous estimate implies that if $v$ belongs to $L^r(t,T;L^s(\Omega))$, not only $ u $ is in $ L^r(t,T;L^s(\Omega))$, but it is also possible to estimate $ \|u\|_{L^r(t,T;L^s(\Omega))}$ by a positive constant independent of the initial datum $u_0$.

Finally, by Theorem \ref{asinto} we know that $u-v$ belongs to  $L^{\overrightarrow{p}}\left(  t,T;W_{0}^{1,\overrightarrow
{p}}(\Omega)\right) $, for every $t \in (0,T)$, and consequently if
 $v$ belongs to $L^{\overrightarrow{p}}\left(  t,T;W_{0}^{1,\overrightarrow
{p}}(\Omega)\right) $ then it follows that $u \in L^{\overrightarrow{p}}\left(  t,T;W_{0}^{1,\overrightarrow
{p}}(\Omega)\right)$ and thus \eqref{regdu} holds true.

\qed



\section*{Acknowledgments}
The authors are partially supported by GNAMPA of the Italian INdAM (National Institute of High Mathematics). The research of G. di Blasio has been funded under  the project "Start" within the program of the University "Luigi Vanvitelli" reserved to young researchers, Piano strategico 2021-2023, moreover it is partially supported by the "Geometric-Analytic Methods for PDEs and Applications" project - funded by European Union - Next Generation EU  within the PRIN 2022 program (D.D. 104 - 02/02/2022 Ministero dell'Universit\`{a} e della Ricerca). This manuscript reflects only the authors' views and opinions and the Ministry cannot be considered responsible for them.


\begin{thebibliography}{999999}


\bibitem{antontsev-chipot-08} S. Antontsev, M. Chipot, \emph{Anisotropic
equations: uniqueness and existence results}, Diff. Int. Eq.
\textbf{21} (2008), 401--419.

\bibitem{ADS} S. N. Antontsev, J. I. D\'{\i}az, S. Shmarev,
\emph{Energy methods for free boundary problems. Applications to nonlinear
PDEs and fluid mechanics}. Progress in Nonlinear Differential Equations and
their Applications, 48. Birkh\"{a}user Boston, Inc., Boston, 2002.




\bibitem{AdBF} A. Alberico ,G. di Blasio, F. Feo,
\emph{Comparison results for nonlinear anisotropic parabolic problems}, Atti Accad. Naz. Lincei Rend. Lincei Mat. Appl. , (2017), 305–322.


\bibitem{AFTL} A. Alvino, V. Ferone, G. Trombetti, P. L. Lions,
\emph{Convex symmetrization and applications}, Ann. Inst. H.
Poincar\'{e} Anal. Non Lin\'{e}aire 14 (1997), 275--293.

\bibitem{AP} D. G. Aronson, L. A. Peletier,
\emph{ Large time behaviour of solutions of the porous medium equation in bounded domains}, J. Diff. Eqns. {\bf 39} (1981), 378-412.

\bibitem{BFK} M. Belloni, V. Ferone, B. Kawohl, \emph{Isoperimetric
inequalities, Wulff shape and related questions for strongly
nonlinear elliptic equations,} Zeit. Angew. Math. Phys. (ZAMP)
54 (2003), 771-789.

\bibitem{BK0} M. Bendahmane, K. Karlsen, \emph{Nonlinear anisotropic elliptic and parabolic equations in  $R^N$  with advection and lower order terms and locally integrable data}\textit{.} Potential Anal. 22 (2005), 207–227.

\bibitem{BK} M. Bendahmane, K. Karlsen, \emph{Anisotropic doubly nonlinear
degenerate parabolic equations}\textit{.} Numerical mathematics and advanced
applications, 381--386, Springer, Berlin, 2006.

\bibitem{BCP} Ph. Benilan, M.G. Crandall, M. Pierre, \emph{ Solutions of the porous medium  in $\R^N$ under optimal conditions on initial values},
             Indiana Univ. Math. V. 33 (1984), 51-87.


\bibitem{BG89} L. Boccardo, T. Gallouet,  \emph{Nonlinear elliptic and parabolic equations involving measure data}, J. Functional Anal.   87  (1989), 149-169.

\bibitem{BMS} L. Boccardo, P. Marcellini, C. Sbordone,\emph{ }$L^{\infty}
$\emph{-regularity for variational problems with sharp nonstandard
growth conditions}, Boll. Un. Mat. Ital. A 4 (1990),
219-225.

\bibitem{BG} M. Bonforte, G. Grillo, \emph{ Super and ultracontractive bounds for doubly nonlinear evolution equations }, Rev. Mat. Iberoamericana 22 (2006), n.1, 11-129.


\bibitem{BaCri} B. Brandolini, F. C. C\^{\i}rstea,  \emph{Singular anisotropic elliptic equations with gradient-dependent lower order terms}, Nonlinear Differential Equations and Applications,  30 (2023).

\bibitem{Barbara} B. Brandolini, F. C. C\^{\i}rstea,  \emph{Boundedness of solutions to singular anisotropic elliptic equations}, Discrete Contin. Dyn. Syst. Ser. S (2023). Doi: 10.3934/dcdss.2023190.

\bibitem{cianchi immersione} A. Cianchi,\emph{ A fully anisotropic Sobolev
inequality}\textit{,} Pacific J. Math. 196 (2000),
283--295.

\bibitem{cianchi anisotropo} A. Cianchi, \emph{Symmetrization in
anisotropic elliptic problems}\textit{,} Comm. Part. Diff. Eq.
\textbf{32} (2007), 693--717.



\bibitem{CG} F. Cipriani, G. Grillo, \emph{ Uniform bounds for solutions to quasilinear parabolic equations}, J. Differential Equations {\bf 177} (2001), 209-234.

\bibitem{DiBH1} E. Di Benedetto, M. A. Herrero,  \emph{ On the Cauchy problem and initial traces for a degenerate parabolic equation}, Trans. AMS 314 (1989), 187-224. 

\bibitem{DiBH2} E. Di Benedetto, M. A. Herrero,  \emph{ Non negative solutions of the evolution p-Laplacian equation. Initial traces and Cauchy problem when $1<p<2$}, Arch. Rational Mech. Anal. 111, 3, (1990), 225-290.

\bibitem{dBFZ}   G. di Blasio, F. Feo, G. Zecca, \emph{Regularity results for local solutions to some anisotropic elliptic equations}, Isr. J. Math. (2023), https://doi.org/10.1007/s11856-023-2564-y.

\bibitem{dBFZ_2}   G. di Blasio, F. Feo, G. Zecca, \emph{Existence and uniqueness of solutions to some anisotropic  elliptic equations with a singular convection term}, in press on Journal d'Analyse Mathématique.

\bibitem{dBL} G. di Blasio, P.D. Lamberti, \emph{Eigenvalues of the Finsler p-Laplacian on varying domains}, Mathematika,
66 (2020), 765-776.


\bibitem{castro}A. Di Castro, \emph{Existence and regularity results for
anisotropic elliptic problems}\textit{,} Adv. Nonlinear Stud.
9 (2009), 367--393.

\bibitem{DMV}  F.G. D\"{u}zg\"{u}n, S. Mosconi, V. Vespri,
\emph{Anisotropic Sobolev embeddings and the speed of propagation for parabolic equations}, J. Evol. Equ. (2019), 845–882.

\bibitem{ELM} L. Esposito, F. Leonetti, G. Mingione, \emph{Sharp
regularity for functionals with }$(p,q)$\emph{ growth}, J.Diff.
Equat. 204 (2004), 5-55.

\bibitem{FVV}  F. Feo,  J. L. Vázquez, B. Volzone, \emph{Anisotropic  p -Laplacian evolution of fast diffusion type}, Adv. Nonlinear Stud. (2021), 523–555.

\bibitem{FGK}  I. Fragal\`{a}, F. Gazzola, B. Kawohl, \emph{Existence and nonexistence results
for anisotropic quasilinear elliptic equations}, Ann. I. H. Poincar\'{e} 21 (2004), 715-734.

\bibitem{FGL}I. Fragal\`{a}, F. Gazzola, G. Liebermann, \emph{Regularity
and nonexistence results for anisotropic quasilinear elliptic
equations in convex domains}\textit{,} Disc. Cont. Dynam. Syst.
(2005), 280-286.

\bibitem{Gi} M. Giaquinta, \emph{Growth conditions and regularity, a
counterexample}, Manus. Math. \textbf{59} (1987), 245-248.

\bibitem{HV} M. A. Herrero, J. L. Vazquez, {\em Asymptotic behaviour of the solutions of a strongly nonlinear parabolic problem}, Ann. Fac. Sci., Toulose Math. (5) {\bf 3}, n.2 (1981), 113-127.

\bibitem{LI} F. Leonetti, A. Innamorati, \emph{Global integrability for weak solutions to some anisotropic
elliptic equations}, Nonlinear Analysis 113 (2015), 430-434.

\bibitem{Li} F. Li, \emph{Anisotropic parabolic equations with measure data. II}, Math. Nachr.,(2006), 1585–1596.

\bibitem{LiZ} F. Li, H. Zhao, \emph{Anisotropic parabolic equations with measure data}, J. Partial Differential Equations, (2001), 21–30.

\bibitem{Mar}P. Marcellini, \emph{Regularity of minimizers of integrals of
the calculus of variations with non standard growth conditions},
Arch. Rat. Mech. Anal. 105 (1989), 267-284.

\bibitem{MRSC} A. Mercaldo, J. D. Rossi, S. Segura de Le{\'o}n and C. Trombetti, \emph{Anisotropic} $p, q$ \emph{-Laplacian equations when } $p$ \emph{ goes to } $1$, Nonlinear Anal. 73 (2010), 3546--3560.

\bibitem{M} F. Mokhtari, \emph{Nonlinear anisotropic parabolic equations in } $L^m$, Arab J. Math. Sci. 20 (2014), 1–10.

\bibitem{P3}  M.M. Porzio,  \emph{ $L^\infty$-regularity for degenerate and singular anisotropic
  parabolic equations},  Boll. U.M.I., 11-A (1997), 697-707.

\bibitem{decay} M.M. Porzio, \emph{ On decay estimates},  Journal of Evolution Equations 9 (2009), 561-591.



\bibitem{fortmon} M.M. Porzio, \emph{ Asymptotic behavior and regularity properties of strongly nonlinear parabolic equations},  Annali di Matematica Pura ed Applicata  198 (2019), 1803-1833.

\bibitem{lower} M.M. Porzio, \emph{  On the influence of some absorption terms on the  solutions of nonlinear parabolic equations}, submitted.

\bibitem{debmon} M. M. Porzio, \emph{ Regularity and time behavior of the solutions to weak monotone parabolic equations},  Journal of Evolution Equations  21 (2021) n. 4 3849-3889.


\bibitem{Str} B. Stroffolini, \emph{Global boundedness of solutions of anisotropic variational problems},  Bollettino UMI 7 (1991), 345-352.

\bibitem{tro} M. Troisi, \emph{ Teoremi di inclusione per spazi di Sobolev non isotropi}, Ricerche Mat., 18 (1969), 3-24.

\bibitem{V1} J. L. Vazquez, \emph{ Smoothing and decay estimates for nonlinear diffusion equations}, Oxford University press 2006.

\end{thebibliography}
\end{document}